%% file: main.tex
\documentclass[10pt]{article}

\input{SetUpForArticle}	

\pgfplotsset{compat=1.17}

\begin{document}
\title{Computational issues by interpolating with inverse multiquadrics: a solution}
\author{Stefano De Marchi$^1$%
    \thanks{Electronic address: \texttt{stefano.demarchi@unipd.it}}, Nadaniela Egidi$^2$%
	\thanks{Electronic address: \texttt{nadaniela.egidi@unicam.it}}, Josephin Giacomini$^2$%
	\thanks{Electronic address: \texttt{josephin.giacomini@unicam.it}}, Pierluigi Maponi$^2$%
	\thanks{Electronic address: \texttt{pierluigi.maponi@unicam.it}}, Alessia Perticarini$^2$%
	\thanks{Electronic address: \texttt{alessia.perticarini@unicam.it}; Corresponding author}}
\affil{$^1$ Department of Mathematics “Tullio Levi-Civita”, University of Padua, via Trieste 63, 35121 Padova, Italy}
\affil{$^2$ School of Science and Technology - Mathematics Division, University of Camerino, via Madonna delle Carceri 9, 62032 Camerino, Italy}
\date{}

\maketitle

\begin{abstract}
We consider the interpolation problem with the inverse multiquadric radial basis function. The problem usually produces a large dense linear system that has to be solved by iterative methods. The efficiency of such methods is strictly related to the computational cost of the multiplication between the coefficient matrix and the vectors computed by the solver at each iteration. We propose an efficient technique for the calculation of the product of the coefficient matrix and a generic vector. This computation is mainly based on the well-known spectral decomposition in spherical coordinates of the Green's function of the Laplacian operator. We also show the efficiency of the proposed method through numerical simulations.
\end{abstract}
\section{Introduction}
Radial basis functions (RBFs) are efficient tools in approximation of functions and data, in the solution of many engineering problems~\cite{hardy1990theory}, including applications in image processing~\cite{flusser1992adaptive,carr1997surface}. They are a suitable tool for scattered data interpolation problems~\cite{lazzaro2002radial}, solution of differential equations~\cite{franke1998solving}, machine learning techniques~\cite{poggio1990networks} and other applications~\cite{hardy1990theory}. In particular, the solution of scattered data interpolation problems leads to a linear system that has a dense matrix of coefficients. We can identify three main difficulties in solving the corresponding linear system: (i) high computational cost when a large number of interpolation points is considered, (ii) numerical instability and (iii) a clever choice of the shape parameter.

In the present paper, we deal with the computational cost of the solution of the linear systems in the case of large interpolation problems. For these linear systems, iterative methods are the standard solvers and the efficiency of these methods can be improved in several ways, for instance, by improving the efficiency of the iterate and/or by providing an effective preconditioner~\cite{Beatson1999} or decreasing the computational cost of each iteration.

The {\it fast multipole method}~\cite{greengard1987fast} has been proposed to reduce the computational cost {{in the solution algorithms for}} integral equations arising in classical scattering theory. This gave rise to a set of effective techniques centred around the principle of \textit{divide et impera}. Among the many existing variants of this idea, several schemes have been developed to provide efficient solution strategies for structured linear systems~\cite{xi2014superfast,xia2012superfast}, eigenvalue problems~\cite{benner2013preconditioned,xi2014fast},{ interpolation problems by RBFs~\cite{cherrie2002fast,Beatson1998,Beatson1992}}, integral equations~\cite{egidi2009efficient,martinsson2005fast,gillman2012direct}, and partial differential equations~\cite{le2006h}. A review of similar techniques can be found in~\cite{cai2018smash}.

In this work, we propose a method aimed to decrease the computational cost of the product of the coefficient matrix and a generic vector. The method is based on a local low-rank representation of the interpolation matrix $A$. This representation resembles ones used in the fast multipole method and is conceptually based on a simple remark: the computational cost of the product of a matrix $A$ of order $N$ and a generic vector is proportional to $N^2$, but if a matrix $A$ admits a decomposition of the form $A=UV$, where $U$ and $V$ are rectangular matrices $N \times p$ and $p \times N$, respectively, the action of $A$ has a computational cost proportional to $2\, N\, p$, so, this decomposition can be profitably used when $p\ll N$. The efficiency of the proposed technique is also shown through numerical experiments.

The paper is organised as follows. In Section~\ref{sec:translation}, we describe the proposed method, in particular, we define the interpolation matrix decomposition, which is based on the spectral expansion of the fundamental solution of the Laplacian operator, and we introduce the translation strategy for the efficient computation of the matrix action. In Section~\ref{sec:simulations}, we show the results of some
numerical experiments. Finally, we give some conclusions and future developments in Section~\ref{sec:conclusions}.

\section{The interpolation problem with inverse multiquadric RBFs}\label{sec:translation}
We consider the interpolation problem where the interpolation function is expressed in terms of inverse multiquadric RBFs. In the following, we focus on iterpolation problem on $\bbr^2$ but the results can be easily generalized in $\bbr ^s$, for different $s$. Initially, we summarise some preliminary concepts and fix notations. Let $\Omega$ be a subset of $\bbr ^2$, $\mathcal{X}=\{\bsx_1,\bsx_2, \dots, \bsx_N\} \subset \Omega$ be the set of $N$ distinct points, usually called \textit{data sites}, and $f_i \in \bbr$ be the \textit{data values}. Moreover, data values are supposed to be obtained from some unknown function $f:\Omega \rightarrow \bbr$ evaluated at the data sites, i.e. 
\begin{equation} \label{interpolation_cond}
f(\bsx_i)=f_i, \qquad i=1, 2, \dots, N.
\end{equation}
Given a RBF defined through $\phi :[0,\infty )\rightarrow \bbr$, $\phi(r_j)=\phi(\Vert\bx-\bsy_j\Vert)$ for $r_j=\Vert\bsx-\bsy_j\Vert$, the interpolation problem consists in finding the approximation function $P_f(\bsx)$ as follows
\begin{equation}\label{interp_eq}
    P_f(\bsx)=\sum_{j=1}^N c_j\phi(||\bsx-\bsy_j||), \qquad \bsx \in \Omega,
\end{equation}
where $||\cdot ||$ is the Euclidean norm, $\bsy_j$, $j=1, \dots, N$, are the centres of the RBFs and coincide with the data sites $\bsy_j=\bsx_j$, and $c_j$, $j=1, \dots, N$, are unknown coefficients. We want to determine these unknowns in such a way that:
\begin{equation}\label{P_fiEqfi}
    P_f(\bsx_i)=f_i,\qquad i=1,2,\ldots,N.
\end{equation}
In order to compute $P_f$ in~\eqref{interp_eq}, we have to solve the system of linear equations~\eqref{P_fiEqfi} which has the form:
\begin{equation}\label{interp_system}
A\bc=\bbf, 
\end{equation}
where the interpolation  matrix $A$ has entries  $a_{ij}=\phi(||\bsx_i-\bsy_j||)$, $i,j=1, \dots, N$, $\bc=(c_1, \dots, c_N)^T$ is the unknown coefficient vector, and $\bbf=(f_1,\dots,f_N)^T$ is the known term given by the interpolation condition~\eqref{P_fiEqfi}.

The linear system~\eqref{interp_system} has a dense coefficients matrix, moreover, for large values of $N$ direct methods usually cannot be used to calculate its numerical solution due to the memory and CPU-time resources required. Thus, iterative methods remain the unique numerical alternative for dealing with such kind of linear systems. In general, the efficiency of such methods depends on the spectral and sparsity properties of the coefficient matrix. The former usually influence the number of iterations required to achieve a given accuracy in the numerical solution, while the latter influence the computational cost of each iteration, which is strictly dependent on the iterative method considered but, for large dense linear systems, most of this cost is due to the computation of the product of the corresponding coefficient matrix and the generic tentative solution. We point out that, from standard arguments on numerical linear algebra, the computational cost of this matrix-vector product is proportional to $N^2$, when $N$ is the order of the matrix, so we need to supply a more efficient technique for this computation.

For the description of the proposed method, we consider the interpolation problem in $\rtwo$ with inverse multiquadric (IMQ-) RBFs, even if, this study can be extended to other families of RBFs and in other dimensions. More precisely, let $\bsy \in \rtwo$ be the center of the IMQ-RBFs, $\bsx \in \rtwo$, and let $t\in \bbr $ be the shape parameter, we consider:
\begin{equation}\label{IMQ_def}
 \phi(||\bsx-\bsy||)=\dfrac{1}{\sqrt{t^2+||\bsx-\bsy||^2}}.
\end{equation}
In the following sections, we introduce the proposed technique. In particular, in Section~\ref{sec:decomposition}, we present the spectral expansion in spherical coordinates of the Green's function of the Laplacian operator and its use in the local decomposition of the matrix $A$. In Section~\ref{sec:translationTec}, we discuss the translation strategy for the efficient computation of the product of the matrix $A$ and a generic vector. Then, in Section~\ref{sec:costoComp}, we analyse the computational cost of the proposed technique. 

\subsection{The Green's function of the Laplacian operator }\label{sec:decomposition}
\begin{figure}[!hbt]
    \centering
     \subfloat[\label{subf:err_trend_a}]{
    \includegraphics[width=.48\textwidth]{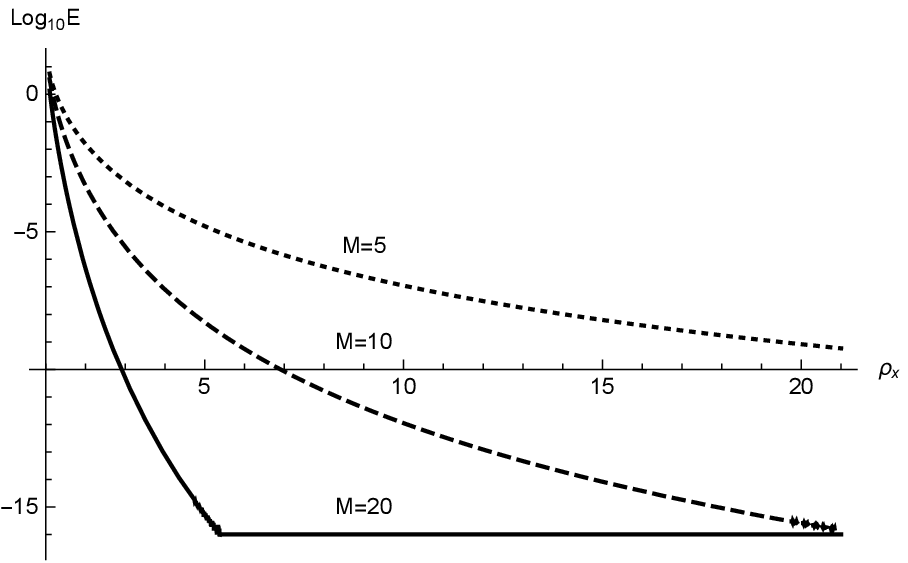}}\hfill
     \subfloat[\label{subf:err_trend_b}]{
     \includegraphics[width=.48\textwidth]{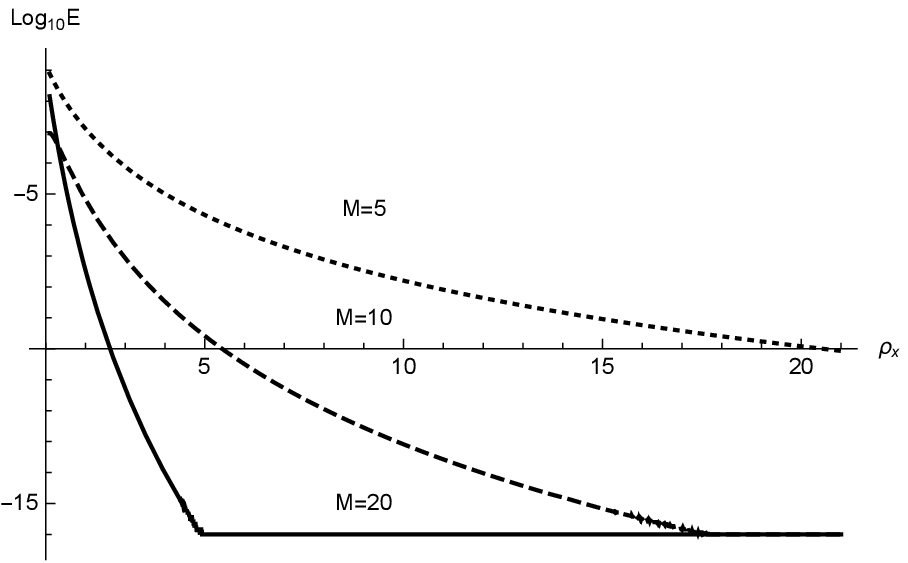}
     }
    \caption{The error $E(\bsX,\bsY)$ in the representation formula~\eqref{serieTroncata} with truncation index $M = 5$, $M=10$, $M = 20$, and $\rho_y=1$, when $\theta_x=\theta_y=\omega_x=\omega_y=\dfrac{\pi}{3}$~\protect\subref{subf:err_trend_a}, $\theta_x=\pi$, $\theta_y=\dfrac{\pi}{4}$, $\omega_x=\dfrac{\pi}{3}$, $\omega_y=\dfrac{\pi}{2}$~\protect\subref{subf:err_trend_b}. We note that the logarithmic scale is used for the axis of the ordinates.}
    \label{fig:ErrorTrend}
\end{figure}

The Green's function, $G(\bsx;\bsy)$, of a linear differential operator $L$ is a solution of the equation $ L G(\bsx;\bsy)=\delta(\bsx-\bsy)$, where $\delta$ is the Dirac delta. From standard arguments, we have that
\begin{equation}\label{G_def}
   G(\bsX;\bsY)=\frac1{|| \bsX-\bsY||}, \qquad \bsX, \bsY \in \bbr^3, \bsX\not=\bsY,
\end{equation}
is the Green's function of the Laplacian operator $L=\Delta$~\cite{carslaw1959conduction}, 
moreover, $G(\bsX;\bsY)$ has a well-known spectral decomposition in spherical coordinates. In particular if $(\rho_x,\theta_x,\omega_x),$  $(\rho_y,\theta_y,\omega_y)$ denote the spherical coordinates of $\bsX$ and $\bsY$, respectively, where $\rho_x,\rho_y \in [0,+\infty), \theta_x,\theta_y \in [0,\pi], \omega_x,\omega_y \in [0,2\pi)$, we have the following spectral expansion of $G(\bsX;\bsY)$ when $\rho_x\ne \rho_y$
\begin{equation}\label{GreenDecompositionG}
\begin{aligned}
   G(\bsX;\bsY)= &\sum_{n=0}^\infty \sum_{m=0}^n \epsilon_m \frac{(n-m)!}{(n+m)!}P_n^m\left(\cos\theta_y\right) P_n^m\left(\cos\theta_x\right) \cdot \\
    & \cdot \cos\left(m(\omega_x-\omega_y)\right)
    \begin{cases}
        \rho_x^n/\rho_y^{n+1}, & \text{if } \rho_y>\rho_x, \\
        \rho_y^n/\rho_x^{n+1}, & \text{if } \rho_x>\rho_y,
    \end{cases}
\end{aligned}
\end{equation}
where $\epsilon_m$ is the Neumann factor, that is $\epsilon_0=1,\epsilon_m=2,\,m>0$, and $P_n^m$ is the associated Legendre function of order $ m $ and degree $ n $, see~\cite{morse1953methods} for more details. 
Furthermore, in a practical application of~\eqref{GreenDecompositionG}, we have to consider the truncated series, that is
\begin{equation}\label{serieTroncata}
\begin{aligned}
  G(\bsX;\bsY)\approx &\sum_{n=0}^M \sum_{m=0}^n \epsilon_m \frac{(n-m)!}{(n+m)!}P_n^m\left(\cos\theta_y\right) P_n^m\left(\cos\theta_x\right) \cdot \\
    & \cdot \cos\left(m(\omega_x-\omega_y)\right)
    \begin{cases}
        \rho_x^n/\rho_y^{n+1}, & \text{if } \rho_y>\rho_x, \\
        \rho_y^n/\rho_x^{n+1}, & \text{if } \rho_x>\rho_y,
    \end{cases}
\end{aligned}
\end{equation}
where $M\in \bbn$ is the truncation parameter for the index $n$ in~\eqref{GreenDecompositionG}. The accuracy of approximation~\eqref{serieTroncata} strongly depends on the particular choice of $\bsX$, $\bsY \in \bbr ^3$. In particular, it depends on the ratio between $\rho_y$ and $\rho_x$, 
as stated in the following theorem.
\begin{theorem}\label{thm:series_error}
Let $\bsX,\bsY\in \bbr^3$ be such that $\rho_y<\rho_x$ and let $r=\dfrac{\rho_y}{\rho_x}$. For the absolute error $E$ in the truncated series~\eqref{serieTroncata}, that is
\begin{equation*}\label{ErrDecomp}
    E(\bsX,\bsY)=\Bigg|G(\bsX;\bsY) - \sum_{n=0}^M \sum_{m=0}^n \epsilon_m \frac{(n-m)!}{(n+m)!}P_n^m\left(\cos\theta_y\right) P_n^m\left(\cos\theta_x\right)\cos\left(m(\omega_x-\omega_y)\right)\dfrac{1}{\rho_x}r^n\Bigg| ,
\end{equation*}
we have the following bound
\begin{equation}\label{maggiorazioneErr}
 E(\bsX,\bsY)\leq \dfrac{1}{\rho_x}\dfrac{r^{M+1}}{(1-r)}.
\end{equation}
\begin{proof}
Let $\alpha$ be the angle between $\bsX$ and $\bsY$, we have:
\begin{equation}
  \begin{aligned}
  E(\bsX,\bsY)=&\Bigg| \sum_{n=M+1}^\infty \sum_{m=0}^n \epsilon_m \frac{(n-m)!}{(n+m)!}P_n^m\left(\cos\theta_y\right) P_n^m\left(\cos\theta_x\right)\cos\left(m(\omega_x-\omega_y)\right)\dfrac{1}{\rho_x}r^n\Bigg|\\
  &= \Bigg| \dfrac{1}{\rho_x}\sum_{n=M+1}^\infty P_n\left(\cos\alpha\right)r^n \Bigg|\leq  \dfrac{1}{\rho_x}\sum_{n=M+1}^\infty\Bigg| P_n\left(\cos\alpha\right)\Bigg|r^n \leq  \dfrac{1}{\rho_x}\sum_{n=M+1}^\infty r^n =\\
  &=\dfrac{1}{\rho_x}\Bigg(\sum_{n=0}^\infty r^n-\sum_{n=0}^M r^n \Bigg) = \dfrac{1}{\rho_x}\Bigg(\dfrac{1}{1-r}-\dfrac{1-r^{M+1}}{1-r}\Bigg) = \dfrac{1}{\rho_x}\dfrac{r^{M+1}}{1-r}
  \end{aligned}  
\end{equation}
where $P_n$ is the Legendre polynomial of degree $n$ and the second equality holds because
\begin{equation*}
    \displaystyle{P_n(\cos \alpha)=\sum_{m=0}^n \epsilon_m \frac{(n-m)!}{(n+m)!}P_n^m\left(\cos\theta_y\right) P_n^m\left(\cos\theta_x\right)\cos\left(m(\omega_x-\omega_y)\right)},
\end{equation*}
see~\cite[Formula (10.3.38)]{morse1953methods} for more details, and the last inequality holds because $|P_n(x)|\leq 1$ for all $\vert x \vert \leq 1$~\cite[Formula (28.35)]{spiegelschaum}.
\end{proof}
\end{theorem}
We note that a similar theorem holds in the case $\rho_x<\rho_y$ and $r=\dfrac{\rho_x}{\rho_y}$. As a consequence of  Theorem~\ref{thm:series_error}, we have that the series in~\eqref{GreenDecompositionG} converges uniformly on every compact set of the domain $||\bsX||>||\bsY||$ or on the domain $||\bsX||<||\bsY||$. In Figure~\ref{fig:ErrorTrend}, we report the result of a numerical experiment showing that, when $\bsX$ is sufficiently far from $\bsY$, formula~\eqref{serieTroncata} is accurate even with moderate values of $M$, while it produces large errors when $||\bsX||\approx ||\bsY||$, even if a large truncation index $M$ is used. More precisely, in Figure~\ref{fig:ErrorTrend}, $\rho_x\in \bigg [\dfrac{11}{10},21\bigg ]$
and $\rho_y=1$; in addition, $\theta_x=\theta_y=\omega_x=\omega_y=\dfrac{\pi}{3}$ in Figure~\ref{fig:ErrorTrend}\subref{subf:err_trend_a}, whereas $\theta_x=\pi$, $\theta_y=\dfrac{\pi}{4}$, $\omega_x=\dfrac{\pi}{3}$ and $\omega_y=\dfrac{\pi}{2}$ in Figure~\ref{fig:ErrorTrend}\subref{subf:err_trend_b}.

\subsection{The translation technique}\label{sec:translationTec}

If we consider $\bsX=(x_1,x_2,t/2)$, $\bsY=(y_1,y_2,-t/2)$, where $t$ is the shape parameter in~\eqref{IMQ_def}, $\bsx=(x_1,x_2)\in\bbr^2$, and $\bsy=(y_1,y_2)\in\bbr^2$, we have that
\begin{equation}\label{G_eq_phi}
G(\bsX;\bsY)=\frac1{|| \bsX- \bsY||}=\frac1{\sqrt{||\bsx-\bsy ||^2+t^2}}=\phi(||\bsx-\bsy ||)
\end{equation} 
and Formula~\eqref{GreenDecompositionG} gives the spectral decomposition for the IMQ-RBFs with shape parameter $t$. We note that
Formula~\eqref{GreenDecompositionG} resembles the usual representation of degenerate kernels, in fact, it expresses $\phi(||{\bsx}-{\bsy}||)$ as a series, where each term is a product of two functions of only one variable, i.e., a function of $\bsx$ and the other of $\bsy$. In addition, it identifies two regions where we have two different representations. In fact, the structure of~\eqref{GreenDecompositionG} requires us to divide the elements of $\bA$ into three sets: the elements where $||\bsx||<||\bsy||$, those where $||\bsx||>||\bsy||$, and those where $||\bsx|| = ||\bsy||$ for which~\eqref{GreenDecompositionG} cannot be used. This is an interesting formula resembling a low-rank decomposition of matrix $A$, but unfortunately it does not provide a low-rank representation of the matrix $A$ due to the aforementioned partition of $A$. However, formula~\eqref{GreenDecompositionG} provides a local decomposition for $\phi(||\bsx - \bsy||)$ that can be profitably used in the solution of~\eqref{interp_system}.

When the ratio $\rho_y/\rho_x \approx 1$ (but $\bsX$ is substantially different from $\bsY$), a simple translation operation can be used in formula~\eqref{serieTroncata}. 
Let $\bsz=(z_1,z_2)\in \Omega \subset \bbr^2$ be such that $||\bsy-\bsz||<||\bsx-\bsz||$, and let $\bsZ=(z_1,z_2,-t/2)$, then from~\eqref{serieTroncata} we have:
\begin{equation}\label{dec_phi_trasl}
\begin{aligned}
   &\phi(||{\bsx}-{\bsy}||)= \phi(||({\bsx}-{\bsz})-({\bsy}-{\bsz})||)=
   G(\bsX-\bsZ;\bsY-\bsZ)\approx\\
   &\approx\sum_{n=0}^M \sum_{m=0}^n \epsilon_m \frac{(n-m)!}{(n+m)!}P_n^m\left(\cos(\theta_{y-z})\right) P_n^m\left(\cos(\theta_{x-z})\right)  \cos\left(m(\omega_{x-z}-\omega_{y-z})\right)  \dfrac{\rho_{y-z}^n}{\rho_{x-z}^{n+1}},
   \end{aligned}
\end{equation}
where $(\rho_{x-z},\theta_{x-z},\omega_{x-z}), \ (\rho_{y-z},\theta_{y-z},\omega_{y-z})$ are the spherical coordinates of vectors $\bsX-\bsZ$, $\bsY-\bsZ$ $\in \bbr^3$, respectively.
The translation technique for the computation of $A\bc$ consists in a convenient use of~\eqref{dec_phi_trasl}.
For simplicity of notation, let us define the following quantities: $d_{n,m}=\epsilon_m \dfrac{(n-m)!}{(n+m)!}$, $h_{n,m}(\bsX-\bsZ)=\dfrac{P_n^m\left(\cos(\theta_{x-z})\right)}{\rho_{x-z}^{n+1}} $, $j_{n,m}(\bsY-\bsZ)=P_n^m\left(\cos(\theta_{y-z})\right)\rho_{y-z}^n$. Thus, when $||\bsy-\bsz||<||\bsx-\bsz||$, we can separate the contribution of ${\bsx}$ and ${\bsy}$ in the following way:
\begin{equation}\label{dec_separata}
\begin{aligned}
   \phi(||{\bsx}-{\bsy}||)=& \sum_{n=0}^M \sum_{m=0}^n d_{n,m}h_{n,m}(\bsX-\bsZ)\cos(m\omega_{x-z})j_{n,m}(\bsY-\bsZ)\cos(m\omega_{y-z})+\\+&\sum_{n=0}^M \sum_{m=0}^n d_{n,m}h_{n,m}(\bsX-\bsZ)\sin(m\omega_{x-z})j_{n,m}(\bsY-\bsZ)\sin(m\omega_{y-z}).
   \end{aligned}
\end{equation}

\begin{figure}[!hbt]
     \subfloat[\label{subf:liv1}]{ \includegraphics[width=0.45\textwidth]{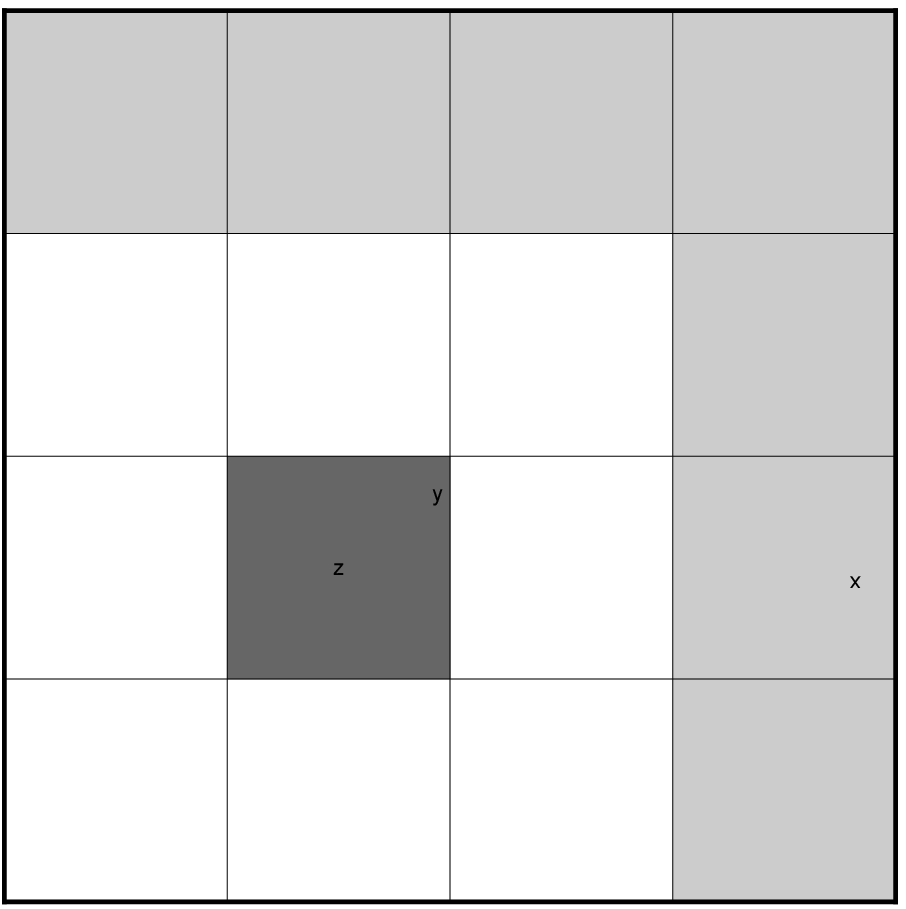}
     }\hfill
    \subfloat[\label{subf:liv2}]{ \includegraphics[width=0.45\textwidth]{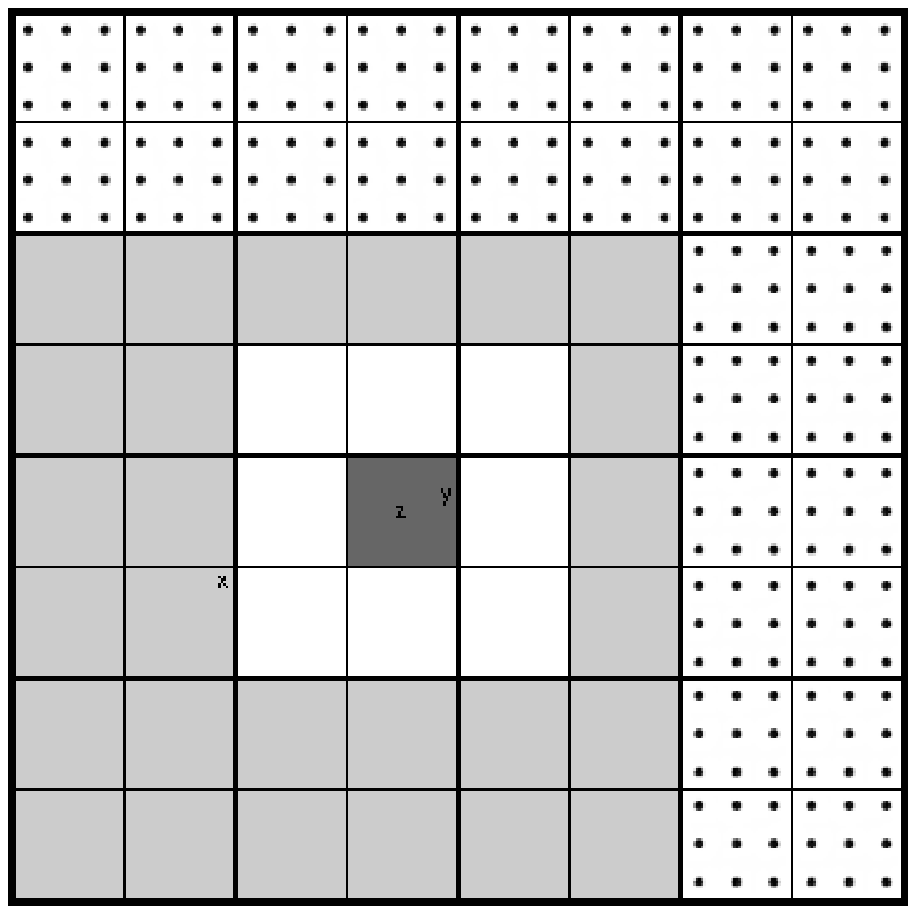}}
    \caption{Example of the translation technique at levels $l=1$~\protect\subref{subf:liv1} and $l=2$~\protect\subref{subf:liv2}.}
    \label{fig:trasla}
\end{figure}

We can suppose, without losing generality, that the domain $\Omega$ is contained in a square  $D\subset \rtwo$. The proposed strategy considers a set of partitions of $D$: the first partition is obtained by dividing $D$ into $4\times 4$ equivalent blocks, the successive are obtained by bisecting the edges of the previous blocks. More precisely, the partition considered in $D$ depends on a recursive index $l$, with $l=1,2,\ldots,L$. For each level $l$, the square D is partitioned in $2^{l+1}\times2^{l+1}$ blocks like in Figure~\ref{fig:trasla}, where the first two levels of these partitions are considered. The idea of the proposed strategy is the following: when ${\bsy}\in \mathcal{X}$ is in a set $S$ (the dark gray squares in Figure~\ref{fig:trasla}) of the considered partition, we associate to ${\bsy}$ the center ${\bsz}$ of $S$ and we use formula~\eqref{dec_phi_trasl} only for ${\bsx}\in \mathcal{X}$ belonging to well-separated squares of the considered partition (the light gray squares in Figure~\ref{fig:trasla}). In this way, we obtain a high accuracy of formula~\eqref{dec_phi_trasl}, even for small values of the truncation index $M$.

In more detail, for each point ${\bsy}$ in a block at level $1$, we choose ${\bsz}$ as the center of this block as shown in the dark gray block in Figure~\ref{fig:trasla}\subref{subf:liv1}. Then, we define the well-separated blocks, with respect to the {{ previously selected block that contains $\bsy$}}, as those blocks of points $\bsx$ such that $\dfrac{\rho_{y-z}}{\rho_{x-z}}<\dfrac{1}{2}$. In Figure~\ref{fig:trasla}\subref{subf:liv1}, the well-separated blocks are the light gray ones. So, {{for each ${\bsy}$,}} in a given block of the partition we apply formula~\eqref{dec_phi_trasl} for points ${\bsx}$ in the corresponding well-separated blocks of the same partition. In this way,{{ if $t=1$ and $D$ has unitary edge, we have an error $E<2.86 \cdot 10^{-9}$ when $M= 10$, and an error $E<4.41 \cdot 10^{-17}$ when $M= 20$. In fact, it is easy to verify that $\rho_{y-z}\leq \dfrac{\sqrt{2}}{8}$ and $\rho_{x-z}\geq\dfrac{\sqrt{73}}{8}$, so that $\dfrac{\rho_{y-z}}{\rho_{x-z}}\leq 0.1655$.}} At the second level, i.e. $l=2$, each block of level $1$ is partitioned in $2\times2$ equivalent squares, thus obtaining $64$ blocks, as in Figure~\ref{fig:trasla}\subref{subf:liv2}. At level $2$, the same procedure is applied to the part of the domain $D$ that has not already been considered at level $1$. So, for each point ${\bsy}$ in a smaller block, as the dark gray block in Figure~\ref{fig:trasla}\subref{subf:liv2}, we choose ${\bsz}$ as the center of the block, and we apply formula~\eqref{dec_phi_trasl} for each point ${\bsx}$ in the new well-separated blocks, i.e., the light gray blocks in Figure~\ref{fig:trasla}\subref{subf:liv2}. {{Also at this level we have the same accuracy of previous level.}} The dotted region in Figure~\ref{fig:trasla}\subref{subf:liv2} 
highlights the region where formula~\eqref{dec_phi_trasl} has been already applied and so it does not need to be treated at this level. We repeat this procedure for all the blocks until the finest level $L$ is reached, where the contributions of the remaining white regions are directly computed by using formula~\eqref{IMQ_def} in the matrix action. 

In Algorithm~\ref{algo}, we report the translation technique for computing the{{ matrix-vector product
$\bb=A\bu$, when $\bu \in \bbr^N$ is a generic vector.
In this algorithm, we use the following notation: $P_l$ is the set of all blocks at level $l$, $p$ is one of these blocks, 
$S_p$ is the set of the blocks well-separated from the block $p$, 
and ${NS}_p$ is the set of the blocks which are not well-separated from the block $p$.

}}
\begin{algorithm}[!hbt]
 \SetAlgoLined
 $\bb=\left(0,\dots,0\right)^T$\;
\For{ $l=1,\dots,L$}{
    \For{ $p\in P_l$}{
        ${\bsz}$ is the center of $p$\\
        \For{${\bsy}_j \in p$}{ 
            $\bv=\left(0,\dots,0\right)^T$\;
            $\bw=\left(0,\dots,0\right)^T$\;
            Compute the contribution of ${\bsy}_j$ from formula~\eqref{dec_separata}\\
            \For{$n=0,\dots,M$}{
                \For{$m=0,\dots,n$}{
                    $k=\dfrac{n(n+1)}{2}+m+1$\\
                    $\bv_k=\bv_k + j_{n,m}(\bsY_j-\bsZ)\sin(m\omega_{y_j-z})\bu_j$\\
                    $\bw_k=\bw_k + j_{n,m}(\bsY_j-\bsZ)\cos(m\omega_{y_j-z})\bu_j$
                }
            }
        }
        \For{${\bsx}_i \in S_p$}{
            Compute the contribution of ${\bsx}_i$ from formula~\eqref{dec_separata}\\
            \For{$n=0,\dots,M$}{
                \For{$m=0,\dots,n$}{
                    $k=\dfrac{n(n+1)}{2}+m+1$\\
                    $\bb_i=\bb_i + d_{n,m} h_{n,m}(\bsX_i-\bsZ)\big[\sin(m\omega_{x_i-z})\bv_k +\cos(m\omega_{x_i-z})\bw_k\big]$
                }
            }
        }
    }
}
\For{$p\in P_L$}{
     \For{${\bsx}_i \in NS_p$}{
         \For{${\bsy}_j \in p$}{
        Compute the contribution of ${\bsx}_i$ and ${\bsy}_j$ from formula~\eqref{interp_system}\\
    $\bb_i=\bb_i + A_{i,j}\bu_j$
        }
    }
}
 \caption{Given $\bu \in \bbr ^N$, computes $\bb=A\bu \in \bbr ^N$ as follows.}\label{algo}
\end{algorithm}

\subsection{The computational cost}\label{sec:costoComp}

The recursive procedure described in the previous section gives an efficient method for computing $A\bu$, where $\bu\in\bbr^N$ is a generic vector and the data sites{{ in $\mathcal{X}$ are generic points}} of the domain $\Omega$. However, for the sake of simplicity, we suppose a uniform distribution of points, so at the generic level $l$, each block $p$ contains $\dfrac{N}{4^{l+1}}$ data sites. As illustrated in Algorithm~\ref{algo}, we have that
\begin{equation}\label{MatxVec}
    \begin{aligned}
    \sum_{j=1}^N A_{i,j}\bu_j=& \scaleto{\sum_{l=1}^L \sum_{\substack{p\in P_l\\i\in S_p}}}{35pt} \left\{ \sum_{n=0}^M \sum_{m=0}^n \Bigl(d_{n,m}h_{n,m}(\bsX_i-\bsZ)\cos(m\omega_{x_i-z})\Bigr) \left( \sum_{j\in p} j_{n,m}(\bsY_j-\bsZ)\cos(m\omega_{y_j-z})\bu_j\right)+\right.\\
    &\phantom{\{}+\left.\sum_{n=0}^M \sum_{m=0}^n\left( d_{n,m}h_{n,m}(\bsX_i-\bsZ)\sin(m\omega_{x_i-z})\right)\left( \sum_{j\in p} j_{n,m}(\bsY_j-\bsZ)\sin(m\omega_{y_j-z})\bu_j\right)\right\}+\\
    &+\sum_{\substack{p\in P_L\\i\in {NS}_p}} \left\{\sum_{j\in p} A_{i,j}\bu_j\right\}, \qquad  \qquad i=1,\dots,N,
   \end{aligned}
\end{equation}
where $i\in S_p$ denotes that $\bsx_i$ is in well-separated blocks $S_p$, while $j\in p$ denotes that $\bsy_j$ is in $p$. This shows that the computation consists of two main addenda. The first addendum in~\eqref{MatxVec}, i.e., the one for $i\in S_p$, $p\in P_l$, $l=1,\dots,L$, can be computed efficiently since the part depending on the row indices $i$ is independent of the part depending on the column indices $j$; while the second addendum has to be evaluated only at the last level $L$. The computational cost $c$ of the matrix-vector product in~\eqref{MatxVec} can be obtained by summing the cost of the various steps. We retrace the computations in Algorithm~\ref{algo} and we evaluate the overall computational cost by counting only the multiplication operations. The computation of the coefficients $d_{n,m}$ and functions $h_{n,m}, j_{n,m}$ as well as trigonometric functions in formula~\eqref{MatxVec} are not considered since they are computed only once at the beginning of the solution process and stored in arrays.

Let $c_1$ be the cost of lines $5-16$ in Algorithm~\ref{algo}, so $c_1$ is equal to
the cost of 
\begin{equation}\label{cost:part1}
   \bv_{n,m}(\bsZ)=\sum_{j\in p} j_{n,m}(\bsY_j-\bsZ)\sin(m\omega_{y_j-z})\bu_j \quad \text{ and} \quad \bw_{n,m}(\bsZ)=\sum_{j\in p} j_{n,m}(\bsY_j-\bsZ)\cos(m\omega_{y_j-z})\bu_j,
\end{equation}
for all $n=0,1,\dots,M,\ m=0,1,\dots,n$, that is $4KN_p$, where $N_p$ are the number of points of $\mathcal{X}$ contained in $p$ and  $K=(M+1)(M+2)/2$.
Thus, considering a block $p$ at a generic level $l$, we have
\begin{equation*}
    c_1\leq 4\frac{N}{4^{l+1}}K,
\end{equation*}

Let $c_2$ be the cost of lines $2-27$ in Algorithm~\ref{algo}, that is the cost of
\begin{equation}\label{cost:part2}
    \sum_{n=0}^M \sum_{m=0}^n d_{n,m}h_{n,m}(\bsX_i-\bsZ)\sin(m\omega_{x_i-z})\bv_{n,m}(\bsZ) \quad \text{and} \quad \sum_{n=0}^M \sum_{m=0}^n d_{n,m}h_{n,m}(\bsX_i-\bsZ)\cos(m\omega_{x_i-z})\bw_{n,m}(\bsZ),
\end{equation}
for all $n=0,1,\dots,M,\ m=0,1,\dots,n$, $p\in P_l,l=1,2,\dots,L,$ $\bsz$ the centre of $p$, $i\in S_p$, where $\bv_{n,m},\bw_{n,m}$ contain the contributions of the points $\bsy_j$ for the sine and cosine part, respectively, already calculated in~\eqref{cost:part1}; $\bv_{n,m},\bw_{n,m}$ correspond respectively to $\bv_k,\bw_k$ in Algorithm~\ref{algo}. In the following discussion, the cost of the multiplication by the factor $d_{n,m}$ is neglected since such factor can be included in $h_{n,m}$ or, equivalently, in $j_{n,m}$ during the construction of the data structures. We analyse the cost $c_2$ as function of the level $l$. At the first level, we have $N/16$ points of $\mathcal{X}$ in each block $p\in P_1$ and at most there are $12$ well-separated blocks from $p$, that are blocks in $S_p$. 
Thus, the cost $c_2^{(1)}$ at the first level is
\begin{equation*}
    c_2^{(1)}\leq 4K 12 \frac{N}{16}=3KN,
\end{equation*}
where $K$ is defined above and gives the number of addenda in each sum appearing in~\eqref{cost:part2}. 
At level $l\geq2$, each block of the level $l-1$ is divided into 4 blocks, so fixing the block $p$ of the first level, $p$ is divided into $4^{l-1}$ blocks at level $l$. 
For each small block there are at most $27$ well-separated blocks. Thus, the cost $c_2^{(l)}$ at level $l$ is
\begin{equation*}
    c_2^{(l)}\leq 4K 4^{l-1} 27 \frac{N}{4^{l+1}}=\frac{27}{4}KN.
\end{equation*}
Since in the costs $c_2^{(l)},l=1,2,\dots,L,$ we fixed a block of the first level and we referred to it during the calculations, by considering also that there are $16$ different blocks at the first level, 
we have
\begin{equation*}
    c_2\leq 16\sum_{l=1}^L c_2^{(l)}\leq 108 KNL.
\end{equation*}

Let $c_3$ be the cost of lines $28-35$ in Algorithm~\ref{algo}, that is the cost of
\begin{equation*}
    \sum_{j\in p} A_{i,j}\bu_j,
\end{equation*}
for all $i\in {NS}_p$ and $p\in P_L$. With analogous arguments to the ones used above, i.e., the number of blocks in the first level, the number of blocks obtained at level $L$ from the subdivision of a first-level block into smaller blocks and the number of points into a block, considering also that the maximum number of non well-separated blocks is at most $9$, we obtain
\begin{equation*}
    c_3\leq 16 \, 4^{L-1} \frac{N}{4^{L+1}}9\frac{N}{4^{L+1}}=9\frac{N^2}{4^{L+1}}.
\end{equation*}
So, the total cost of the procedure is
\begin{equation}\label{total_cost}
    \begin{aligned}
    c=&c_1+c_2+c_3\leq\\
    \leq&N\left(\frac{4K}{4^{l+1}}+108KL+\frac{9N}{4^{L+1}}\right)\leq\\
    \leq&N\left(109KL+\frac{9N}{4^{L+1}}\right).
    \end{aligned}
\end{equation}
In formula~\eqref{total_cost}, the parameters $K$ and $N$ are fixed before the procedure, whereas the number of levels $L$ should be chosen in order to obtain the minimum $c$. Now, if we consider the upper bound for $c$, we should search for the minimum of this upper bound. In more detail, let $\alpha=109K, \beta=9N/4, g(\lambda)=\displaystyle{\alpha {{\lambda}}+\frac{\beta}{4^\lambda}}$, where $\lambda$ is the real extension of the integer variable $L$. Then, it can be shown that $g$ has a minimum in $\lambda_{min}=\displaystyle{\log_4\left(\frac{\beta\log_{10}(4)}{\alpha}\right)}$, where $g(\lambda_{min})=\displaystyle{\frac{\alpha}{\log_{10}(4)}\left(\log_{10}\left(\frac{\beta\log_{10}(4)}{\alpha}\right)+1\right)}$. A similar result can be obtained without the extension of the integer variable $L$, because of the convexity of the function $g$. In conclusion, the upper bound of $c$ is proportional to $N\log_{10}(N)$.

\section{Numerical simulations}\label{sec:simulations}
We present a numerical experiment performed with the technique described in the previous section, in order to show the accuracy and the efficiency of the proposed method. The experiment consists in computing the product of the matrix $A$ and a random vector whose components range is $[-1,1]$, by using two different techniques: the standard matrix-vector multiplication, and the translation technique described in Section~\ref{sec:translation}.

We fix the shape parameter $t=1$, and we consider the set of quasi-random $N$ Halton points~\cite{HALTON1960} in $\Omega=[0,1]^2$, with $N=2(4), 4(4), 6(4), 8(4), 1(5)$, where $x(y)$ denotes the real number $x \cdot 10^y$. Furthermore, in the translation technique, the truncation index in~\eqref{dec_phi_trasl} is fixed to $M = 10$ and the maximum number of levels is $L = 1,2,3$. The results are reported in Table~\ref{tab:results} and Figure~\ref{fig:tempi}. Table~\ref{tab:results} shows the elapsed time $T$ in the computation of the product of $A$ and a random vector by using a standard row-column product, the elapsed time in the same computation by using the translation technique, the relative error $E$ in infinity-norm between the result computed by the translation technique and the one obtained by the standard row-column product, and the number of levels $L$ that gives the best efficiency result among $L=1,2,3$. In addition, Figure~\ref{fig:tempi} reports a comparison of the execution times of the different methods. In the $x-$axis the number of interpolation points is reported, while the vertical axis gives the normalised execution time, that is the execution time divided by the number of interpolation points, i.e., $\dfrac{T}{N}$. The line with circular markers represents the usual matrix-vector multiplication, while the line with squared markers represents the proposed technique. In Figure~\ref{fig:tempi}, the substantial reduction of the execution time can be easily appreciated, as well as the logarithmic trend of the proposed strategy. 

The results in Table~\ref{tab:results} and Figure~\ref{fig:tempi} have been obtained on a Workstation equipped with an Intel(R) Xeon(R) CPU E5-2620 v3 @2.40GHz, operative system Red Hat Enterprise Linux, release 7.5. All computations have been made in double precision and the FORTRAN codes have been compiled by the NAGWare f95 Compiler.
\begin{table}[!hb]
    \centering
    \begin{tabular*}{450pt}{@{\extracolsep\fill}lcccc@{\extracolsep\fill}}
        \toprule
        \multirow{2}{*}{$N$} & \text{Standard} &  \multicolumn{3}{c}{\text{Translation}} \\
        \cmidrule(lr){2-2}\cmidrule(lr){3-5}
        & $T$ \, [\text{\SI{}{\second}}] & $T$ \, [\text{\SI{}{\second}}] & $E$ & $L$ \\
        \midrule 
        \vspace{2pt}
        $2(4)$ & $1.62(0)$ & $1.58(0)$ & $2.67(-9)$ & $1$ \\
        \vspace{2pt}
        $4(4)$ & $1.79 (1)$ & $5.23(0)$ & $4.61(-9)$ & $2$ \\
        \vspace{2pt}
        $6(4)$ & $3.21 (1)$ & $9.78(0)$ & $6.62(-9)$ & $2$ \\
        \vspace{2pt}
        $8(4)$ & $8.38(1)$ & $1.58(1)$ & $8.72(-9)$ & $2$ \\
        \vspace{2pt}
        $1(5)$ & $1.28(2)$ & $2.33(1)$ & $1.06(-8)$ & $2$ \\
       
        \bottomrule
        \end{tabular*}
        \caption{Execution times $T$ (in seconds), relative error $E$ and number of levels $L$ by varying the number of points $N$, for the truncation index $M=10$.\label{tab:results}}
    \end{table}
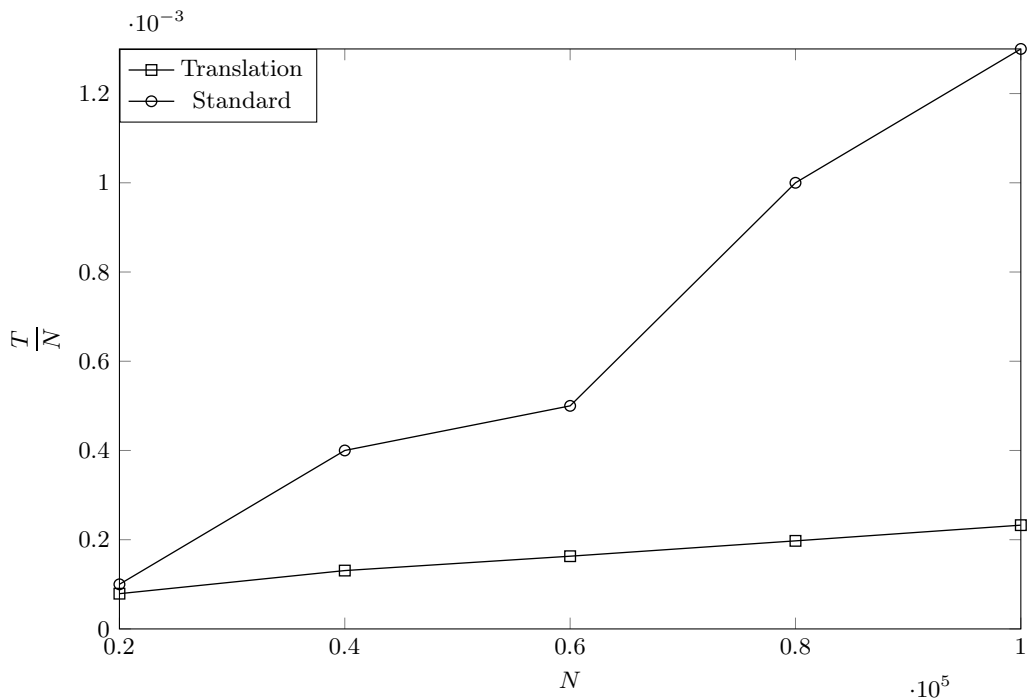
\begin{figure}[!hbt]
\centering
\begin{tikzpicture}
\begin{axis} [scale only axis,height=.3\textheight,width=.7\textwidth,
    xmin=20000,xmax=100000,ymin=0,ymax=0.0013,xlabel={$N$},ylabel={$\dfrac{T}{N}$},xtick={20000,40000,60000,80000,100000},legend style={at={(0,1)},anchor=north west}]

\addplot+ [color=black, mark=square, line width=.5pt] table[x index=0, y index=1] {tempiEsecuzione/fast.txt};

\addlegendentry{Translation}
\addplot+ [color=black, mark=o, line width=.5pt] table[x index=0, y index=1] {tempiEsecuzione/classico.txt};
\addlegendentry{Standard}
\end{axis}
\end{tikzpicture}
 \caption{Trends of the normalised execution times with the two strategies, when $M=10$}
    \label{fig:tempi}
\end{figure}

From these results, we can observe a substantially reduced computational time $T$ of the translation technique with respect to the standard computation, especially for large value of $N$. Moreover, the error $E$ shows the high accuracy provided by the proposed method despite the small value of the truncation index $M$. 

\section{Conclusions}\label{sec:conclusions}
We considered a scattered interpolation problem with the IMQ-RBFs and we proposed an efficient strategy for computing the product between the coefficient matrix and a generic vector. This method is based on the well-known decomposition formula in spherical coordinates for the Laplacian operator and a simple translation technique. The presented numerical experiments show strongly promising results, both for the efficiency and the accuracy of the proposed technique. In fact, such computational strategy has a smaller computational cost than the standard matrix-vector multiplication, and the computed numerical solutions are almost the same. These results encourage further investigations about this method. In particular, such strategy should be implemented into an iterative method for the solution of the interpolation problem. Then, it could be interesting to study the generalization of this technique to the interpolation problem with different choices of RBFs, as well as to the collocation problem for the solution of differential equations. Finally, the decomposition of $A$ could be profitably exploited for the construction of ad-hoc preconditioner for the interpolation problem.

\vskip 0.2in
{\bf Acknowledgments} This research has been accomplished within Rete ITaliana di Approssimazione (RITA), the thematic group on "Approximation Theory and Applications" of the Italian Mathematical Union and partially funded by GNCS-IN$\delta$AM.
\printbibliography
\end{document}

%% file: SetUpForArticle.tex
\usepackage[T1]{fontenc}                  

\usepackage[utf8]{inputenc}     

\usepackage{microtype}                     

\usepackage[a4paper, margin=0.8in]{geometry}	

\usepackage[italian,english]{babel}	

\usepackage[babel]{csquotes}
                                           
\usepackage{emptypage}			

\usepackage{indentfirst}			

\usepackage{booktabs}			

\usepackage[dvipsnames]{xcolor}	

\usepackage{tabularx}			
\usepackage{multirow}

\usepackage{subfig}				

\usepackage{graphicx}			
\graphicspath{{./Figures/}}

\usepackage[export]{adjustbox}

\usepackage{pgfplots}
\pgfplotsset{legend pos=south east, every axis/.append style={font=\small}}

\usepackage{caption}			

\usepackage{listings}			

\usepackage[font=small]{quoting}	

\usepackage{amsmath,amssymb,amsthm,mathrsfs,scalerel}        

\usepackage[euler]{textgreek}		

\usepackage{mathtools}			

\usepackage[linesnumbered,ruled,vlined]{algorithm2e}
\SetKw{Break}{break}      

\usepackage{siunitx}  

\usepackage[english]{varioref}		

\usepackage{mparhack,relsize}		

\usepackage[style=numeric-comp,citestyle=ieee,hyperref,backref,natbib,backend=biber]{biblatex} 

\bibliography{Bibliography}

\usepackage[affil-it]{authblk}		

\usepackage{eurosym}			

\usepackage[bookmarks=true]{hyperref}		

\usepackage{bookmark}			

\usepackage{empheq}			

\usepackage{lineno}				

\usepackage{comment}

\newcommand{\bu}{\mathbf{u}}
\newcommand{\bx}{\mathbf{x}}

\newcommand{\bb}{\mathbf{b}}
\newcommand{\bc}{\mathbf{c}}
\newcommand{\bv}{\mathbf{v}}
\newcommand{\bw}{\mathbf{w}}

\newcommand{\bbf}{\mathbf{f}}

\newcommand{\bA}{\mathbf{A}}

\newcommand{\rtwo}{\mathbb{R}^2}
\newcommand{\bbr}{\mathbb{R}}
\newcommand{\bbn}{\mathbb{N}}

\newcommand{\bsx}{\boldsymbol{x}}
\newcommand{\bsy}{\boldsymbol{y}}
\newcommand{\bsz}{\boldsymbol{z}}
\newcommand{\bsX}{\boldsymbol{X}}
\newcommand{\bsY}{\boldsymbol{Y}}
\newcommand{\bsZ}{\boldsymbol{Z}}

\providecommand{\keywords}[1]{\textit{Keywords ---} #1}

\theoremstyle{definition}

\theoremstyle{plain}
\newtheorem{theorem}{Theorem}[section]

\theoremstyle{remark}

\newcommand\restr[2]{{
  \left.\kern-\nulldelimiterspace 
  #1 
  \right|_{#2} 
  }}



%% file: main.bib
@book{morse1953methods,
  title={Methods of Theoretical Physics. Vol. 1-2},
  author={Morse, Philip M and Feshbach, Herman},
  year={1953},
  publisher={New York},
  pages={1274}
}

@article{HALTON1960,
author = {HALTON, J.H.},
journal = {Numerische Mathematik},
pages = {84-90},
title = {On the efficiency of certain quasi-random sequences of points in evaluating multi-dimensional integrals.},
url = {http://eudml.org/doc/131448},
volume = {2},
year = {1960},
}

@article{carr1997surface,
  title={Surface interpolation with radial basis functions for medical imaging},
  author={Carr, Jonathan C and Fright, W Richard and Beatson, Richard K},
  journal={IEEE transactions on medical imaging},
  volume={16},
  number={1},
  pages={96--107},
  year={1997},
  publisher={IEEE}
}

@article{lazzaro2002radial,
  title={Radial basis functions for the multivariate interpolation of large scattered data sets},
  author={Lazzaro, Damiana and Montefusco, Laura B},
  journal={Journal of Computational and Applied Mathematics},
  volume={140},
  number={1-2},
  pages={521--536},
  year={2002},
  publisher={Elsevier}
}

@article{franke1998solving,
  title={Solving partial differential equations by collocation using radial basis functions},
  author={Franke, Carsten and Schaback, Robert},
  journal={Applied Mathematics and Computation},
  volume={93},
  number={1},
  pages={73--82},
  year={1998},
  publisher={Elsevier}
}

@article{poggio1990networks,
  title={Networks for approximation and learning},
  author={Poggio, Tomaso and Girosi, Federico},
  journal={Proceedings of the IEEE},
  volume={78},
  number={9},
  pages={1481--1497},
  year={1990},
  publisher={IEEE}
}

@article{Beatson1999,
  title={Fast fitting of radial basis functions: Methods based on preconditioned GMRES iteration},
  author={Beatson, Richard K and Cherrie, Jon B and Mouat, Cameron T},
  journal={Advances in Computational Mathematics},
  volume={11},
  number={2},
  pages={253--270},
  year={1999},
  publisher={Springer}
}

@article{Beatson1998,
  title={Fast evaluation of radial basis functions: Moment-based methods},
  author={Beatson, Richard Keith and Newsam, Garry Neil},
  journal={Siam Journal on Scientific Computing},
  volume={19},
  number={5},
  pages={1428--1449},
  year={1998},
  publisher={SIAM}
}

@article{Beatson1992,
  title={Fast evaluation of radial basis functions: I},
  author={Beatson, Richard K and Newsam, Garry N},
  journal={Computers \& Mathematics with Applications},
  volume={24},
  number={12},
  pages={7--19},
  year={1992},
  publisher={Elsevier}
}

@book{spiegelschaum,
  title={Schaum's Outlines: Mathematical Handbook of Formulas and Tables},
  author={Spiegel, Murray R and Lipschutz, Seymour and Liu, John},
  %volume={2},
  pages={166},
  year={2218},
  publisher={McGraw-Hill}
}

@techreport{carslaw1959conduction,
  title={Conduction of heat in solids},
  author={Carslaw, Horatio Scott and Jaeger, John Conrad},
  year={1959},
  institution={Clarendon press,}
}

@article{hardy1990theory,
  title={Theory and applications of the multiquadric-biharmonic method 20 years of discovery 1968--1988},
  author={Hardy, Rolland L},
  journal={Computers \& Mathematics with Applications},
  volume={19},
  number={8-9},
  pages={163--208},
  year={1990},
  publisher={Elsevier}
}

@article{flusser1992adaptive,
  title={An adaptive method for image registration},
  author={Flusser, Jan},
  journal={Pattern Recognition},
  volume={25},
  number={1},
  pages={45--54},
  year={1992},
  publisher={Elsevier}
}

@article{greengard1987fast,
  title={A fast algorithm for particle simulations},
  author={Greengard, Leslie and Rokhlin, Vladimir},
  journal={Journal of computational physics},
  volume={73},
  number={2},
  pages={325--348},
  year={1987},
  publisher={Elsevier}
}

@article{cai2018smash,
  title={SMASH: Structured matrix approximation by separation and hierarchy},
  author={Cai, Difeng and Chow, Edmond and Erlandson, Lucas and Saad, Yousef and Xi, Yuanzhe},
  journal={Numerical Linear Algebra with Applications},
  volume={25},
  number={6},
  pages={e2204},
  year={2018},
  publisher={Wiley Online Library}
}

@article{egidi2009efficient,
  title={The efficient solution of direct medium problems by using translation techniques},
  author={Egidi, Nadaniela and Maponi, Pierluigi},
  journal={Mathematics and Computers in Simulation},
  volume={79},
  number={8},
  pages={2361--2372},
  year={2009},
  publisher={Elsevier}
}

@article{cherrie2002fast,
  title={Fast Evaluation of Radial Basis Functions: Methods for Generalized Multiquadrics in $\backslash$RR\^{}$\backslash$protectn},
  author={Cherrie, Jon Barry and Beatson, Richard Keith and Newsam, Garry N},
  journal={SIAM Journal on Scientific Computing},
  volume={23},
  number={5},
  pages={1549--1571},
  year={2002},
  publisher={SIAM}
}

@article{xi2014superfast,
  title={Superfast and stable structured solvers for Toeplitz least squares via randomized sampling},
  author={Xi, Yuanzhe and Xia, Jianlin and Cauley, Stephen and Balakrishnan, Venkataramanan},
  journal={SIAM Journal on Matrix Analysis and Applications},
  volume={35},
  number={1},
  pages={44--72},
  year={2014},
  publisher={SIAM}
}

@article{xia2012superfast,
  title={A superfast structured solver for Toeplitz linear systems via randomized sampling},
  author={Xia, Jianlin and Xi, Yuanzhe and Gu, Ming},
  journal={SIAM Journal on Matrix Analysis and Applications},
  volume={33},
  number={3},
  pages={837--858},
  year={2012},
  publisher={SIAM}
}

@article{benner2013preconditioned,
  title={The preconditioned inverse iteration for hierarchical matrices},
  author={Benner, Peter and Mach, Thomas},
  journal={Numerical Linear Algebra with Applications},
  volume={20},
  number={1},
  pages={150--166},
  year={2013},
  publisher={Wiley Online Library}
}

@article{xi2014fast,
  title={A fast randomized eigensolver with structured LDL factorization update},
  author={Xi, Yuanzhe and Xia, Jianlin and Chan, Raymond},
  journal={SIAM Journal on Matrix Analysis and Applications},
  volume={35},
  number={3},
  pages={974--996},
  year={2014},
  publisher={SIAM}
}

@article{martinsson2005fast,
  title={A fast direct solver for boundary integral equations in two dimensions},
  author={Martinsson, Per-Gunnar and Rokhlin, Vladimir},
  journal={Journal of Computational Physics},
  volume={205},
  number={1},
  pages={1--23},
  year={2005},
  publisher={Elsevier}
}

@article{gillman2012direct,
  title={A direct solver with O (N) complexity for integral equations on one-dimensional domains},
  author={Gillman, Adrianna and Young, Patrick M and Martinsson, Per-Gunnar},
  journal={Frontiers of Mathematics in China},
  volume={7},
  number={2},
  pages={217--247},
  year={2012},
  publisher={Springer}
}

@article{le2006h,
  title={H-matrix preconditioners in convection-dominated problems},
  author={Le Borne, Sabine and Grasedyck, Lars},
  journal={SIAM journal on matrix analysis and applications},
  volume={27},
  number={4},
  pages={1172--1183},
  year={2006},
  publisher={SIAM}
}
